\documentclass[11pt]{amsart}
\usepackage[T1]{fontenc}
\usepackage{lmodern,amsmath,amssymb,amsthm,booktabs,microtype}
\usepackage[margin=1in]{geometry}
\usepackage[colorlinks=true,linkcolor=blue,citecolor=blue,urlcolor=blue]{hyperref}
\newtheorem{theorem}{Theorem}[section]
\newtheorem{proposition}[theorem]{Proposition}
\newtheorem{lemma}[theorem]{Lemma}
\newtheorem{corollary}[theorem]{Corollary}
\theoremstyle{remark}
\newtheorem{remark}[theorem]{Remark}
\newcommand{\Z}{\mathbb Z}
\newcommand{\CP}{\mathbb{CP}^{2}}
\newcommand{\CPb}{\overline{\mathbb{CP}}{}^{2}}

\title[Exotic universal covers and free cyclic actions]{Exotic Universal Covers and Free Cyclic Actions
\newline from the Cartwright-Steger Surface}

\author{Anar Akhmedov}
\address[A.~Akhmedov]{School of Mathematics, University of Minnesota, Minneapolis, MN 55455, USA}
\curraddr{Department of Mathematics, Harvard University, Cambridge, MA 02138, USA}
\email{akhmedov@umn.edu}

\author{B. Doug Park}
\address[B.~D.~Park]{Department of Pure Mathematics, University of Waterloo, Waterloo, ON, N2L 3G1, Canada}
\email{bdpark@uwaterloo.ca}

\author{S\"{u}meyra Sakall{\i}}
\address[S.~Sakall{\i}]{Department of Mathematics and Statistics, University of South Florida, 4202 E. Fowler Avenue, Tampa, FL 33620, USA}
\email{sumeyras@usf.edu}

\author{Sai-Kee Yeung}
\address[S.-K.~Yeung]{Department of Mathematics, Purdue University, West Lafayette, IN 47907, USA}
\email{yeungs@purdue.edu}

\date{October 25, 2026}

\begin{document}
\raggedbottom

\begin{abstract}
For every $m\geq2$, we construct a closed minimal symplectic four-manifold with fundamental group $\Z/m$, Euler characteristic $35$, and signature $1$. The construction is a genus-five symplectic sum of a blowup of the regular fourfold Cartwright-Steger cover and a cyclic Akhmedov-Park manifold. The gluing identifies the homomorphisms from the two surface groups to $\Z/m$. The universal cover is minimal symplectic and is homeomorphic but not diffeomorphic to $(18m-1)\CP\#(17m-1)\CPb$. It admits a free symplectic action of $\Z/m$. Intermediate quotients give examples with fundamental group $\Z/n$ and characteristic numbers $(e,\sigma)=(35k,k)$ for every $n\geq2$ and $k\geq1$.
Using the same Cartwright-Steger summand, we also construct, for every $m\geq2$, a minimal symplectic four-manifold with fundamental group $\Z/m$, Euler characteristic $31$, and signature $1$, whose minimal symplectic universal cover is an exotic copy of $(16m-1)\CP\#(15m-1)\CPb$ with a free symplectic $\Z/m$-action.
\end{abstract}
\maketitle

\section{Introduction}

We construct minimal symplectic four-manifolds with arbitrary nontrivial finite cyclic fundamental group and fixed Euler characteristic and signature. Our first construction uses a nodal curve on the Cartwright-Steger surface and the cyclic manifolds of Akhmedov-Park \cite[Section~9]{AP}. We choose the genus-five gluing so that the homomorphisms from the two surface groups to $\Z/m$ agree.

\begin{theorem}\label{thm:main}
For every integer $m\geq2$, there is a closed minimal symplectic four-manifold $Z_m$ such that
\[
\pi_1(Z_m)\cong\Z/m,\qquad
(e(Z_m),\sigma(Z_m),\chi_h(Z_m),c_1^2(Z_m))=(35,1,9,73).
\]
Its universal cover $\widetilde Z_m$ is minimal symplectic and is homeomorphic, but not diffeomorphic, to
\[
(18m-1)\CP\#(17m-1)\CPb.
\]
The deck group acts freely and symplectically on $\widetilde Z_m$. Both $Z_m$ and $\widetilde Z_m$ have odd intersection form on integral second homology modulo torsion and are nonspin.
\end{theorem}

Here $\chi_h=(e+\sigma)/4$ and $c_1^2=2e+3\sigma$. The notation $\CPb$ denotes the complex projective plane with the opposite orientation. We use the symplectic orientation throughout.

Using complex ball quotients, Akhmedov, Sakall\i\ and Yeung \cite{ASY} constructed simply connected symplectic four-manifolds with $(e,\sigma)=(35,1)$. We use their genus-five surface construction with the cyclic Akhmedov-Park manifold $B_m$. 
The other summand comes from the regular fourfold Albanese cover of the Cartwright-Steger surface.  
Its surface subgroup is determined in \cite[Propositions~6.3--6.4 and Lemma~6.8]{AYcurve}. The choice of gluing retains the group $\Z/m$.

Our earlier paper \cite{APSY} gives cyclic quotients with $(e,\sigma)=(43,1)$ by a genus-two sum followed by a torus-bundle sum. The construction here gives $(e,\sigma)=(35,1)$. General symplectic constructions with prescribed fundamental group are due to Gompf \cite{Gompf}; related results for cyclic groups appear in \cite{BSS,LT}.

For $m=2$, the homeomorphism type of the universal cover already occurs in \cite{ASY}. Our construction gives a free symplectic involution whose quotient has signature one. We do not determine whether the smooth structure on its universal cover agrees with that of any of those earlier examples.

\begin{corollary}\label{cor:main-ray}
For every $n\geq2$ and $k\geq1$, there is a closed minimal symplectic four-manifold $W_{n,k}$ with
\[
\pi_1(W_{n,k})\cong\Z/n,\qquad
(e,\sigma,\chi_h,c_1^2)(W_{n,k})=(35k,k,9k,73k).
\]
Its universal cover is homeomorphic but not diffeomorphic to
$(18nk-1)\CP\#(17nk-1)\CPb$.
\end{corollary}

In Theorem~\ref{thm:smaller-cyclic} we retain the Cartwright-Steger summand and replace $B_m$ by a blowup of $\Sigma_2\times T^2$. We perform four Luttinger surgeries after the genus-five sum. For every $m\geq2$, this gives a closed minimal symplectic four-manifold with fundamental group $\Z/m$ and $(e,\sigma)=(31,1)$. Its minimal symplectic universal cover is homeomorphic but not diffeomorphic to $(16m-1)\CP\#(15m-1)\CPb$.

The paper is organized as follows.
Sections~\ref{sec:negative} and~\ref{sec:positive}
construct the two symplectic summands.
Section~\ref{sec:sum} describes the gluing and computes
the fundamental groups of the resulting manifolds,
while Section~\ref{sec:universal} studies their
universal covers.
Section~\ref{sec:quotients} proves
Corollary~\ref{cor:main-ray} and records the integral
homology and irreducibility of the quotients.
Finally, Section~\ref{sec:smaller-cyclic} presents
a smaller construction with Euler characteristic $31$.

We write $\Sigma_g$ for a closed oriented surface of genus $g$ and $\nu F$ for a tubular neighborhood of an embedded surface $F$. Our commutator convention is $[x,y]=xyx^{-1}y^{-1}$. The symbol $\#_F$ denotes a symplectic sum along $F$ with the specified normal-bundle gluing.

\subsection*{Acknowledgments}
The second author was partially supported by an NSERC Discovery grant 
RGPIN--2024--04401.
The authors acknowledge the use of ChatGPT (OpenAI)
for assistance with revising the exposition, refining
selected arguments, checking computations, and locating
references. The authors take full responsibility for
the mathematical content and references in this paper.

\section{The Cartwright-Steger pair}
\label{sec:negative}

Let $X$ be the Cartwright-Steger surface, let $\Gamma=\pi_1(X)$, and let $\alpha:X\to E$ be its Albanese map. The presentation and abelianization in \cite{CKY} give
\[
H_1(X;\Z)\cong\Z^2,\qquad
\alpha_*:\Gamma\twoheadrightarrow\Z^2,
\qquad (e(X),\sigma(X))=(3,1).
\]
We write
\[
N_2=\ker\bigl(\Gamma\xrightarrow{\alpha_*}\Z^2
                 \longrightarrow(\Z/2)^2\bigr).
\]
The corresponding regular unramified cover $p:X_4\to X$ has degree four, deck group $(\Z/2)^2$, and
\[
\pi_1(X_4)=N_2,\qquad (e(X_4),\sigma(X_4))=(12,4).
\]
It is a complex ball quotient and hence is aspherical and minimal. This is the regular degree-four Albanese base cover considered in \cite{AYmon}.

We use the type A nodal curve $D\subset X$ of \cite{CKY,ASY}. Its normalization has genus four, it has one ordinary positive node, and $D^2=-1$. The geometric subgroup calculation in \cite[Proposition~6.3]{AYcurve} identifies
\begin{equation}\label{eq:nodalinput}
\operatorname{im}\bigl(\pi_1(D)\longrightarrow\Gamma\bigr)=N_2.
\end{equation}
Here $\pi_1(D)$ is the group of the singular curve, including its branch-joining loop. Equation~\eqref{eq:nodalinput} identifies its image subgroup in $\Gamma$. The branch-joining generator and the exact index calculation are given in \cite[Lemma~6.1 and Proposition~6.3]{AYcurve}. The four components of $p^{-1}(D)\subset X_4$ are described in \cite[Proposition~6.4]{AYcurve}. We recall the smoothing argument.

\begin{lemma}\label{lem:negative-surface}
There is an embedded symplectic surface $D_4\subset X_4$ of genus five and square $-1$ for which
\[
\pi_1(D_4)\longrightarrow\pi_1(X_4)=N_2
\]
is surjective.
\end{lemma}

\begin{proof}
The homomorphism $\Gamma\to(\Z/2)^2$ defining the cover is trivial on $\pi_1(D)$ by \eqref{eq:nodalinput}. Consequently the restriction of the cover to $D$ is trivial: $p^{-1}(D)$ is the disjoint union of four degree-one copies of the nodal curve. Choose one, denoted $D'$. The inclusion $D^{\prime}\subset X_4$ induces a homomorphism with image $N_2$, up to the choice of base path. Normality of $N_2$ makes this independent of the chosen component.

Resolve the positive node of $D'$ symplectically. In a complex coordinate ball the replacement is the usual annulus $uv=\varepsilon$ in place of $uv=0$, adjusted near the boundary. It preserves the homology class and hence the square. It adds one handle to the genus-four normalization, so the resulting surface $D_4$ has genus five and square $-1$.

Let $q_1,q_2$ be the two points of the normalization that map to the node. Smoothing removes small disks about $q_1,q_2$ and joins their boundary circles by an annulus. Loops representing the fundamental group of the normalization may be chosen away from these disks and remain on $D_4$. The added handle has two standard generators. The first follows a path in the punctured normalization from one boundary circle to the other and returns through the annulus. Under inclusion in $X_4$, it represents the same element as the loop in $D^{\prime}$ that passes from one branch of the node to the other. The second generator is the core circle of the annulus. It lies in the coordinate ball and is therefore nullhomotopic in $X_4$. These loops generate the fundamental group of $D_4$ together with the normalization loops. Consequently
\[
\operatorname{im}\bigl(\pi_1(D_4)\to\pi_1(X_4)\bigr)
=\operatorname{im}\bigl(\pi_1(D^{\prime})\to\pi_1(X_4)\bigr)=N_2.
\qedhere \]
\end{proof}

\begin{lemma}\label{lem:negative-complement}
Blow up two distinct regular points of $D_4$ and let
\[
A=X_4\#2\CPb,\qquad F_-\subset A
\]
be its proper transform. Then
\[
g(F_-)=5,\quad F_-^2=-3,\quad (e(A),\sigma(A))=(14,2),
\]
and inclusion induces an isomorphism
\[
\pi_1(A\setminus\operatorname{int}\nu F_-)\cong\pi_1(A)=N_2.
\]
The map $\pi_1(F_-)\to N_2$ is surjective.
\end{lemma}
\begin{proof}
Each blowup decreases the square by one and leaves the genus and ambient fundamental group unchanged. Each exceptional sphere meets $F_-$ once. Removing its small intersection disk leaves a disk in the complement whose boundary is a meridian of $F_-$. The meridian is therefore trivial. The kernel of the complement-to-ambient fundamental-group map is normally generated by a meridian, as follows by filling the tubular neighborhood and applying Seifert-Van Kampen theorem. 
Hence this map is an isomorphism. Lemma~\ref{lem:negative-surface} gives the final surjectivity statement.
\end{proof}

\section{The cyclic positive pair}
\label{sec:positive}

Fix an integer $m\geq2$. Let $a_i,b_i$ ($i=1,2$) be standard generators of $\pi_1(\Sigma_2)$, and let $c,d$ generate the torus factor of $\Sigma_2\times T^2$. Let $\alpha_1,\ldots,\alpha_4$ be the standard generators of $\pi_1(T^4)$. Primes indicate the disjoint parallel copies fixed in \cite{AP}.

We use the notation and the genus-two gluing of \cite[Section~9, Equation~(17)]{AP}. The product-side summing surface is $\Sigma_2\times\{\mathrm{pt}\}\subset\Sigma_2\times T^2$. On the other side, $\overline\Sigma_2\subset T^4\#\CPb$ is the embedded square-zero symplectic genus-two surface constructed in \cite[Section~3]{AP}. Its exceptional sphere meets it twice. With the markings of \cite[Equation~(17)]{AP}, the images of the surface generators under inclusion into the two summands are identified by
\[
a_1\mapsto\alpha_1,\qquad b_1\mapsto\alpha_2,\qquad
a_2\mapsto\alpha_3^2,\qquad b_2\mapsto\alpha_4.
\]
The boundary gluing reverses the orientation of each normal-circle fiber of the summing surfaces. We write
\[
 B_m=\widetilde X_1(m,1)
   =Y_1(1/m,1)\mathbin{\#}_{\Sigma_2=\overline\Sigma_2}Z''(1,1).
\]
The four surgeries in the $\Sigma_2\times T^2$ summand are
\begin{align*}
 &(a_1'\times c',a_1',-1),\qquad
   (b_1'\times c'',b_1',-1),\\
 &(a_2'\times c',c',+1/m),\qquad
   (a_2''\times d',d',+1).
\end{align*}
The other two surgeries, in $T^4\#\CPb$, are
\[
(\alpha_2'\times\alpha_3',\alpha_3',-1),\qquad
(\alpha_2''\times\alpha_4',\alpha_4',-1),
\]
with the disjoint pushoffs specified in \cite[Equation~(6)]{AP}. These define $Z''(1,1)$. All six are Luttinger surgeries. The closed manifold
$B_m$ is symplectic and has
\[
 e(B_m)=5,\qquad\sigma(B_m)=-1,\qquad
 \pi_1(B_m)=\langle c\mid c^m=1\rangle.
\]

\begin{lemma}\label{lem:positive}
The manifold $B_m$ is minimal and contains an embedded symplectic surface
$F_+$ with $g(F_+)=5$ and $F_+^2=3$. The homomorphism $\pi_1(F_+)\to\pi_1(B_m)=\Z/m$ induced by inclusion is surjective.
\end{lemma}

\begin{proof}
Before the surgeries, the genus-two sum is minimal by Usher's theorem \cite{U}. The summand $\Sigma_2\times T^2$ is minimal. Asphericity of $T^4$ implies that every exceptional class in $T^4\#\CPb$ is the signed blowup class, which has nonzero intersection with $\overline\Sigma_2$. Thus the other summand is relatively minimal. Neither summand is ruled. Luttinger surgery preserves minimality
\cite[Proposition 3.1]{HL}, so $B_m$ is minimal.

We construct a genus-two surface $S$ by joining the punctured exceptional sphere to two punctured product tori in the symplectic sum. Let $E\subset T^4\#\CPb$ be the exceptional sphere. We choose the surgery neighborhoods away from the blowup ball, so that $E$ remains unchanged. 
Let $x_1,x_2$ be the images of the two points
$E\cap\overline\Sigma_2$ on the product summing surface $\Sigma_2$. By an isotopy of
the gluing, arrange that both $x_1,x_2$ lie in one small disk in $\Sigma_2$ avoiding the curves
$a_1',b_1',a_2',a_2''$, with paths to the standard product basepoint
in the nonsurgery region. This isotopy does not change the induced gluing on
fundamental groups. Choose the product-side surgery neighborhoods so that their projections under $\Sigma_2\times T^2\to\Sigma_2$ also avoid $x_1,x_2$. Thus both product tori
$P_i=\{x_i\}\times T^2$ are disjoint from all four product-side surgery
neighborhoods.

For each $i$, choose a small disk $\Delta_i^-\subset\overline\Sigma_2$ about the corresponding point of $E\cap\overline\Sigma_2$, and a disk $\Delta_i^+\subset\Sigma_2\times\{\mathrm{pt}\}$ about $x_i$. The gluing identifies $\Delta_i^-$ with $\Delta_i^+$. Arrange that $E$ meets $\overline\Sigma_2$, and each $P_i$ meets $\Sigma_2\times\{\mathrm{pt}\}$, symplectically orthogonally. Choose local trivializations of the normal bundles over these disks so that the normal-circle gluing is complex conjugation followed by a fixed rotation. We make these choices by isotoping the normal-bundle gluing over the disks. Remove two intersection disks from $E$ and one from each $P_i$. At the $i$th attachment, join the punctured $E$ to the punctured $P_i$ by the symplectic annulus in the local sum model $zw=\varepsilon$.
The resulting connected symplectic surface $S$ has
\[
 g(S)=2,\qquad S^2=-1.
\]
For the normal Euler number, trivialize the complex normal direction
tangent to the summing surface over each attachment disk. In the local
model $zw=\varepsilon$, the boundary normal frames match with zero
winding. The relative Euler numbers therefore add to
$E^2+P_1^2+P_2^2=-1$. A normal-bundle twist in the genus-two gluing has
zero winding on these contractible attachment circles. This is the
surface used in the genus-five construction of \cite{ASY}.

A parallel product section $\Sigma'_2$ survives the sum and the surgeries when its $T^2$ coordinate is chosen away from the summing point and all surgery curves. It intersects $S$ in two positive transverse points, one in each $P_i$. Their symplectic resolution is the surface $F_+$. Hence
\[
 g(F_+)=2+2+2-1=5,\qquad F_+^2=-1+0+2\cdot 2=3.
\]
Choose a $c$-loop in $P_1$ avoiding the deleted disk and the resolution disks.
Using the paths in the nonsurgery region, it represents a conjugate
of the standard generator $c$ in $\pi_1(B_m)=\Z/m$. Thus the homomorphism $\pi_1(F_+)\to\Z/m$ induced by inclusion is surjective. Changing the coefficient
from $+1$ to $+1/m$ changes only the boundary identification used to reattach $T^2\times D^2$ in that torus surgery;
the surfaces just constructed are disjoint from that neighborhood.
\end{proof}

\section{The symplectic sum and its fundamental group}
\label{sec:sum}

We use the pair $(A,F_-)$ of Lemma~\ref{lem:negative-complement} and write
$q_-:\pi_1(F_-)\twoheadrightarrow N_2$ for its surface homomorphism.
By definition, $N_2=\alpha_*^{-1}(2\Z^2)$. Surjectivity of $\alpha_*:\Gamma\to\Z^2$ therefore gives $\alpha_*(N_2)=2\Z^2$. Hence the homomorphism
\[
 \vartheta:N_2\longrightarrow\Z,\qquad
 \vartheta(h)=\frac{\alpha_{*,1}(h)}2
\]
is well defined and surjective; $\alpha_{*,1}$ denotes the first
coordinate of $\alpha_*$. Let $\vartheta_m$ be its reduction modulo $m$,
and let $r_m:\pi_1(F_+)\to\Z/m$ be the homomorphism induced by inclusion
in $B_m$. The composition $\vartheta\circ q_-$ is surjective and hence represents a primitive class in $H^1(F_-;\Z)$. The homomorphism $r_m$ also has a primitive integral lift. Indeed, the surviving embedded $c$-loop
maps to one under $r_m$ and is therefore nonseparating: a separating loop is nullhomologous and maps to zero. Extend it to a
symplectic basis of $H_1(F_+;\Z)$, assign coefficient one on this
basis element, and lift the remaining coefficients of $r_m$ to integers.

The action of the mapping class group on first homology surjects onto $\operatorname{Sp}(10,\Z)$. Since this group acts transitively on primitive integral covectors, there is an orientation-preserving diffeomorphism $\phi_m:F_-\to F_+$ matching the two integral lifts. In particular,
\[
 r_m\circ(\phi_m)_*=\vartheta_m \circ q_-.
\]
After rescaling the symplectic forms and applying Moser's theorem, we may take $\phi_m$ to preserve the induced area forms. The opposite normal Euler numbers give a fiber-reversing normal-bundle gluing. 
We form the symplectic sum 
\[
 Z_m=A\mathbin{\#}_{\phi_m}B_m.
\]

\begin{lemma}\label{lem:pushout}
Let $q:G\to H$ be a surjective homomorphism, let $\vartheta:H\to C$ be a homomorphism, and set $r=\vartheta\circ q$. Then the canonical homomorphism
\[
j_C:C\longrightarrow H*_G C
\]
is an isomorphism. Under this identification, the canonical homomorphism from $H$ to the pushout is $\vartheta$.
\end{lemma}

\begin{proof}
Let us write $P=H*_G C$ for the pushout, and let $j_H:H\to P$ and $j_C:C\to P$ be the canonical homomorphisms. The defining relation
\[
j_H\circ q=j_C\circ r=j_C\circ\vartheta\circ q
\]
and the surjectivity of $q$ imply
\[
j_H=j_C\circ\vartheta.
\]
Thus $j_H(H)\subseteq j_C(C)$. Since these two images generate $P$, the homomorphism $j_C$ is surjective.

The homomorphisms $\vartheta:H\to C$ and $\operatorname{id}_C:C\to C$ agree after composition with $q$ and $r$, respectively. By the universal property of the pushout, they induce a homomorphism $f:P\to C$ satisfying
\[
f\circ j_H=\vartheta,
\qquad
f\circ j_C=\operatorname{id}_C.
\]
The second equality makes $j_C$ injective, so $j_C$ is an isomorphism with inverse $f$. The first equality gives the asserted identification of the canonical map from $H$. 
\end{proof}

\begin{proposition}\label{thm:cyclicapplication}
For every $m\geq2$, the gluing just chosen gives a minimal symplectic four-manifold $Z_m$ with
\[
 \pi_1(Z_m)=\Z/m,\qquad
 (e(Z_m),\sigma(Z_m),\chi_h(Z_m),c_1^2(Z_m))=(35,1,9,73).
\]
Its intersection form modulo torsion is odd; hence it is not spin.
\end{proposition}

\begin{proof}
Lemma~\ref{lem:negative-complement} shows that the meridian of $F_-$ is trivial. The gluing identifies a meridian of $F_+$ with the inverse of a meridian of $F_-$. Consequently it kills the meridian $\mu_+$ of $F_+$, and the positive-side group becomes
\[
\pi_1(B_m\setminus\operatorname{int}\nu F_+)/\langle\!\langle\mu_+\rangle\!\rangle
\cong\pi_1(B_m)\cong\Z/m.
\]
The common boundary is a circle bundle over $\Sigma_5$. Its group is generated by lifts $a_i,b_i$ of standard surface generators and a central fiber meridian $\mu$, with relation
\[
 \prod_{i=1}^{5}[a_i,b_i]=\mu^{\pm3}.
\]
The exponent has absolute value three because the normal circle bundles of $F_-$ and $F_+$ have Euler numbers $-3$ and $+3$; its sign depends on the boundary and meridian conventions. This exponent is independent of the choice of normal-bundle gluing. After the meridian is killed, the relation becomes the surface relation. The remaining relations give the pushout
\[
 \pi_1(Z_m)=N_2*_{\pi_1(\Sigma_5)}\Z/m.
\]
After identifying the surface groups by $(\phi_m)_*$, the matching condition
is $r_m=\vartheta_m \circ q_-$, and Lemma~\ref{lem:pushout}
proves the group assertion. A different choice of normal-bundle gluing can multiply the images of lifted surface generators by powers of the meridian. The meridian is trivial, so these factors do not change the group.

The symplectic sum formulas give
\[
 e(Z_m)=14+5+4(5-1)=35,\qquad\sigma(Z_m)=2-1=1.
\]
Consequently $\chi_h(Z_m)=9$ and $c_1^2(Z_m)=2(35)+3=73$.

The positive block is minimal. Under the blowdown $A\to X_4$, every sphere is nullhomotopic in $X_4$ by asphericity. An exceptional sphere therefore has class in the kernel lattice $\langle-1\rangle^2$. Its square is $-1$, so its class is $\pm E_1$ or $\pm E_2$, where $E_1,E_2$ are the blowup classes. Such a class has nonzero
intersection with $F_-$. Neither summand is a ruled surface with the
summing surface as a section. Usher's theorem in \cite{U} gives minimality.

An even unimodular intersection form has signature divisible by eight.
The free intersection form of $Z_m$ is unimodular and has signature one,
so it is odd. A spin manifold has even intersection form, so $Z_m$ is not spin.
\end{proof}

\section{The universal covers}
\label{sec:universal}

\begin{lemma}\label{lem:positive-cover}
The universal cyclic cover $\widehat B_m\to B_m$ is minimal symplectic.
\end{lemma}

\begin{proof}
Before the six surgeries, define a homomorphism to $\Z/m$ by sending $c$ to
one and all other product and $T^4$ generators to zero. Its restrictions from the two summands agree on the
summing surface, so it defines an $m$-fold cover of the unsurgered sum.
The product-side cover is the connected cover
$\Sigma_2\times\widehat T^2$, where the $c$-circle is covered with degree $m$. The other
side consists of  $m$ copies of $T^4\#\CPb$. The summing surface lifts to
$m$ disjoint product sections. Successively sum along them. Each
exceptional sphere intersects its summing surface twice. The product
block is not ruled. 
After the first sum 
\[
\bigl(\Sigma_2\times\widehat{T}^2\bigr)
\#_{\Sigma_2=\overline\Sigma_2}
\bigl(T^4\#\overline{\mathbb{CP}}{}^{2}\bigr),
\] 
the resulting manifold is minimal and has $e=5$ and $\sigma=-1$. 
It is therefore not ruled, since every minimal ruled surface has signature zero.
Usher's theorem in \cite{U} applies
successively and proves the minimality of the unsurgered degree-$m$ cover.

The surgery neighborhoods and their inverse images are disjoint from the genus-two summing surfaces. We perform the lifted surgeries after the successive sums, removing these neighborhoods and reattaching copies of $T^2\times D^2$ along the lifted boundary maps. The homomorphism to $\Z/m$ extends over each downstairs surgery because it sends the prescribed surgery slope, the image of the meridian of the reattached $T^2\times D^2$, to zero. The product-side
$c$ surgery has surgery slope $\mu+mc$. On the connected boundary cover,
\[
\widehat\mu\longmapsto\mu,\qquad
\widehat c\longmapsto mc,
\]
so the lifted slope is $\widehat\mu+\widehat c$. The table below lists the lifted surgeries. The coefficients are measured in the pulled-back Lagrangian framings, with the conventions of \cite{AP}.

\bigskip 

\begin{center}
\begin{tabular}{lcc}
\toprule
Base surgery & number of lifts &  coefficient upstairs\\
\midrule
$(a_1'\times c',a_1',-1)$ & one & $-1$\\
$(b_1'\times c'',b_1',-1)$ & one & $-1$\\
$(a_2'\times c',c',+1/m)$ & one & $+1$\\
$(a_2''\times d',d',+1)$ & $m$ & $+1$ on each\\
each surgery in $Z''(1,1)$ & $m$ & $-1$ on each\\
\bottomrule
\end{tabular}
\end{center}

\bigskip

In the first two rows, write $\lambda=a_1'$ or $b_1'$, respectively. The cover has degree $m$ along the other longitude, the $c$ direction, while $\widehat\lambda\mapsto\lambda$ and $\widehat\mu\mapsto\mu$. Thus the surgery slope $\mu-\lambda$ lifts to $\widehat\mu-\widehat\lambda$, and the coefficient remains $-1$. In the last two rows the homomorphism to $\Z/m$ is trivial on each surgery neighborhood; its inverse image consists of $m$ copies mapping diffeomorphically to that neighborhood, so the coefficients are unchanged. All $3m+3$ lifted torus surgeries are Luttinger
surgeries in the pulled-back Lagrangian framings. Ho-Li's minimality
theorem in \cite{HL} therefore proves minimality after the lifted surgeries. The extended homomorphism $\pi_1(B_m)\to\Z/m$ sends the generator $c$ to one and is therefore an isomorphism. Thus this cover is the universal cover $\widehat B_m$.
\end{proof}

\begin{lemma}\label{lem:cyclicparity}
Let $V$ be a closed oriented four-manifold with finite cyclic fundamental
group and odd intersection form on $H_2(V;\Z)/\operatorname{Tor}$.
Its universal cover $\widetilde V$ is nonspin.
\end{lemma}

\begin{proof}
Let us write $G=\pi_1(V)$ and $p:\widetilde V\to V$. The Hopf exact sequence (cf.~\cite[Chapter VII, \S7, Exercise 6]{Brown})
and $H_2(G;\Z)=0$ give a surjection
\[
 H_2(\widetilde V;\Z)_G\longrightarrow H_2(V;\Z).
\]
The equality $H_2(G;\Z)=0$ follows from the periodic resolution for a finite cyclic group. This map is induced by $p_*$. If $\widetilde V$ were spin, naturality
of $w_2$ would give
$\langle w_2(V),\overline{p_*x}\rangle=\langle w_2(\widetilde V),\bar x\rangle=0$
for every $x\in H_2(\widetilde V;\Z)$, where bars denote reduction modulo two. For $a\in H_2(V;\Z)$, Wu's formula gives
\[
\langle w_2(V),\bar a\rangle=Q_V(a,a)\pmod2,
\]
where $Q_V$ is the intersection pairing. Surjectivity of $p_*$ would therefore make every integral self-intersection in $V$ even, a contradiction.
\end{proof}

\begin{proof}[Proof of Theorem~{\rm \ref{thm:main}}]
Fix $m\ge2$. Proposition~\ref{thm:cyclicapplication} establishes the assertions concerning $Z_m$. We prove the remaining assertions about its universal cover.
Let
\[
p_-:\widehat A\longrightarrow A,
\qquad
p_+:\widehat B_m\longrightarrow B_m
\]
be the connected degree-$m$ covers corresponding respectively to
$\ker\vartheta_m<\pi_1(A)=N_2$ and the trivial subgroup of
$\pi_1(B_m)\cong\mathbb Z/m$. Before removing the tubular neighborhoods of the summing surfaces, set
\[
\widehat F_-=p_-^{-1}(F_-),
\qquad
\widehat F_+=p_+^{-1}(F_+).
\]
The homomorphisms $\vartheta_m\circ q_-$ and $r_m$ are surjective, so both inverse images are connected, and each restriction
$\widehat F_\pm\to F_\pm$ is a degree-$m$ covering. Consequently,
\[
2-2g(\widehat F_\pm)=m(2-2\cdot5)=-8m,
\qquad
g(\widehat F_\pm)=4m+1.
\]
Their normal bundles are the pullbacks of the normal bundles of $F_\pm$, and hence
\[
\widehat F_-^{\,2}=-3m,
\qquad
\widehat F_+^{\,2}=3m.
\]

The matching condition
\[
r_m\circ(\varphi_m)_*=\vartheta_m\circ q_-
\]
ensures that the covering homomorphisms agree under the boundary gluing: both kill the fiber meridian, and their induced homomorphisms on the surface groups agree. Thus the boundary gluing lifts to a diffeomorphism $\widehat\Phi_m$ between the boundaries of the lifted tubular neighborhoods. Gluing the covered complements gives a connected degree-$m$ cover of $Z_m$:
\[
\bigl(\widehat A\setminus\operatorname{int}\nu\widehat F_-\bigr)
\cup_{\widehat\Phi_m}
\bigl(\widehat B_m\setminus\operatorname{int}\nu\widehat F_+\bigr).
\]
Since $\pi_1(Z_m)\cong\mathbb Z/m$, this is its universal cover $\widetilde Z_m$. Lifting the symplectic sum construction therefore gives
\[
\widetilde Z_m
=
\widehat A\#_{\widehat F_-=\widehat F_+}\widehat B_m,
\]
with the symplectic form pulled back from $Z_m$.

Let $X_{4m}\to X_4$ be the connected degree-$m$ cover corresponding to $\ker\vartheta_m$. The two blowup points in $X_4$ have $2m$ inverse images, so
\[
\widehat A\cong X_{4m}\#2m\overline{\mathbb{CP}}^{\,2}.
\]
The manifold $X_{4m}$ is aspherical. If $S\subset\widehat A$ is an embedded sphere with $S^2=-1$, its composition with the blowdown map is nullhomotopic. Its homology class therefore lies in the kernel of the blowdown homomorphism:
\[
[S]=\sum_{i=1}^{2m}a_iE_i,
\]
where $E_1,\ldots,E_{2m}$ are the exceptional classes. Since
\[
-1=[S]^2=-\sum_{i=1}^{2m}a_i^2,
\]
we have $[S]=\pm E_i$ for some $i$. Each exceptional sphere meets $\widehat F_-$ once, so
\[
[S]\cdot[\widehat F_-]=\pm1.
\]
Thus $\widehat A\setminus\widehat F_-$ contains no exceptional sphere. 

The other summand $\widehat B_m$ is minimal by Lemma~\ref{lem:positive-cover}. Moreover,
\[
\sigma(\widehat A)=2m,
\qquad
\sigma(\widehat B_m)=-m,
\]
so neither summand admits an $S^2$-bundle structure, since such bundles have signature zero. Usher’s theorem \cite[Theorem 1.1]{U} now implies that $\widetilde Z_m$ is minimal. 

Euler characteristic and signature multiply under finite coverings. Hence
\[
e(\widetilde Z_m)=35m,
\qquad
\sigma(\widetilde Z_m)=m.
\]
Since $\widetilde Z_m$ is simply connected, these identities give
\[
b_2^+(\widetilde Z_m)=18m-1,
\qquad
b_2^-(\widetilde Z_m)=17m-1.
\]
By Proposition~\ref{thm:cyclicapplication}, the intersection form of $Z_m$ modulo torsion is odd, so Lemma~\ref{lem:cyclicparity} implies that $\widetilde Z_m$ is nonspin. Simple connectivity then implies that its intersection form is odd. It is therefore the odd indefinite unimodular form
\[
(18m-1)\langle1\rangle
\oplus
(17m-1)\langle-1\rangle.
\]
Freedman’s classification \cite{F}, together with the vanishing of the Kirby–Siebenmann invariant for smooth manifolds, shows that $\widetilde Z_m$ is homeomorphic to
\[
(18m-1)\mathbb{CP}^{2}
\#(17m-1)\overline{\mathbb{CP}}^{\,2}.
\] 

Finally, $b_2^+(\widetilde Z_m)>1$, so Taubes’s theorem \cite{T} gives a nonzero Seiberg–Witten invariant. The standard connected sum above has vanishing Seiberg–Witten invariants, since it splits as a connected sum of two manifolds with positive $b_2^+$. Thus the homeomorphism cannot be a diffeomorphism. The deck group $\mathbb Z/m$ acts freely and preserves the pulled-back symplectic form, completing the proof. 
\end{proof}

\section{Quotients and integral homology}
\label{sec:quotients}

\begin{proof}[Proof of Corollary~{\rm \ref{cor:main-ray}}]
Fix integers $n\geq2$ and $k\geq1$, and apply Theorem~\ref{thm:main} with $m=nk$. Let $H$ be the unique order-$n$ subgroup of the deck group $\Z/(nk)$ of $\widetilde Z_{nk}\to Z_{nk}$, and set
\[
W_{n,k}=\widetilde Z_{nk}/H.
\]
Since $\widetilde Z_{nk}$ is simply connected, $\pi_1(W_{n,k})\cong H\cong\Z/n$. The resulting quotient
covers $Z_{nk}$ with degree $k$, so its characteristic numbers are
as displayed. Its universal cover is the same minimal symplectic
manifold. An exceptional sphere in the quotient would lift to an
exceptional sphere upstairs, so the quotient is minimal.
\end{proof}

In summary, the characteristic numbers of $Z_m$ are independent of $m$. Those of its universal cover are obtained by multiplication by $m$.

\bigskip

\begin{center}
\begin{tabular}{@{}lrrrr@{}}
\toprule
Manifold & $e$ & $\sigma$ & $\chi_h$ & $c_1^2$\\
\midrule
cyclic quotient $Z_m$ & $35$ & $1$ & $9$ & $73$\\
universal cover $\widetilde Z_m$ & $35m$ & $m$ & $9m$ & $73m$\\
\bottomrule
\end{tabular}
\end{center}

\bigskip

\noindent
Integral duality gives
\[
\begin{aligned}
H_1(Z_m;\Z)&\cong\Z/m,&
H_2(Z_m;\Z)&\cong\Z^{33}\oplus\Z/m,& H_3(Z_m;\Z)&=0.
\end{aligned}
\]
Indeed, the fundamental groups determine $H_1$, the Euler characteristic determines the free rank of $H_2$, and Poincare duality and the universal coefficient theorem identify its torsion with $\operatorname{Ext}(H_1,\Z)$. Also $H_3\cong H^1=0$.

\begin{corollary}
The manifolds $Z_m$, their universal covers, and the intermediate quotients of Corollary~{\rm \ref{cor:main-ray}} are irreducible: every smooth connected-sum decomposition has a homotopy-sphere summand.
\end{corollary}
\begin{proof}
They are minimal symplectic and have finite cyclic or trivial fundamental group. These groups are residually finite, so the conclusion follows from Hamilton-Kotschick \cite[Corollary~1]{HK}.
\end{proof}

\section{A smaller cyclic construction}
\label{sec:smaller-cyclic}

\begin{theorem}\label{thm:smaller-cyclic}
For every integer $m\geq 2$, there is a closed minimal symplectic four-manifold $V_m$ with
\[
 \pi_1(V_m)\cong\Z/m,
 \qquad (e,\sigma,\chi_h,c_1^2)(V_m)=(31,1,8,65).
\]
Its universal cover is minimal symplectic and is homeomorphic, but not
diffeomorphic, to
\[
 (16m-1)\CP\#(15m-1)\CPb.
\]
The deck group acts freely and symplectically on the universal cover.
\end{theorem}

We use the Cartwright-Steger pair $(A,F_-)$ of
Lemma~\ref{lem:negative-complement}. Let
$q_-:\pi_1(F_-)\twoheadrightarrow N_2$ be the map induced by inclusion and let
$\ell=\vartheta\circ q_-$, where $\vartheta$ is the half-Albanese coordinate
from Section~\ref{sec:sum}. 
Since \(\mathbb Z\) is abelian, \(\ell\) factors through \(H_1(F_-;\mathbb Z)\). 
We use the same notation for the resulting homomorphism \(\ell:H_1(F_-;\mathbb Z)\to\mathbb Z\).

\begin{lemma}\label{lem:smaller-n2-homology}
The abelianization of $N_2$ is $\Z^2\oplus(\Z/2)^3$.
\end{lemma}

\begin{proof}
We use Stover's presentation of $\Gamma$ on generators
$x,y,z$ with relators $R_1,\ldots,R_{12}$
\cite[Section~3]{Stover}. The associated presentation
complex is obtained by attaching twelve disks to a
wedge of three oriented circles labeled $x,y,z$,
with attaching words $R_1,\ldots,R_{12}$.
Its fundamental group is $\Gamma$.

In Stover's abelianization coordinates $r,s$, the
generators have images
\[
x\longmapsto rs^{-1},\qquad
y\longmapsto rs^2,\qquad
z\longmapsto r^{-3}s.
\]
The subgroup $N_2$ determines a degree-four covering
of the presentation complex. The covering has four
vertices, all lying over the single vertex downstairs.
We label them by the elements of $(\Z/2)^2$:
\[
0=(0,0),\qquad 1=(1,0),\qquad
2=(0,1),\qquad 3=(1,1).
\]
Addition of these labels means addition in
$(\Z/2)^2$, with $r,s$ corresponding to $1,2$.
Thus the positive lifts of $x,y,z$ starting at a
vertex $v$ end at $v+3,v+1,v+3$, respectively.

Write $(v,a)$ for the positive lift of the $a$-edge
starting at $v$. There are twelve edges, one for
each generator and each initial vertex. The three
edges
\[
(0,x):0\longrightarrow3,\qquad
(0,y):0\longrightarrow1,\qquad
(1,x):1\longrightarrow2
\]
connect all four vertices without a cycle and
therefore form a maximal tree. Collapsing this
tree leaves nine edges, which we order by initial
vertex and then by the labels $x,y,z$.

Each relator $R_i$ has four closed lifts, one
starting at each vertex. After collapsing the
tree, these give forty-eight relations in the
nine remaining generators. Their exponent sums
form a $48\times9$ integer matrix $M$, with
\[
H_1(N_2;\Z)\cong
\operatorname{coker}\bigl(
M^{\mathsf T}:\Z^{48}\longrightarrow\Z^9
\bigr).
\]
We order the rows first by $R_1,\ldots,R_{12}$
and then by the initial vertices $0,1,2,3$.
To compute a row, follow the lifted relator and
omit the tree edges. For each non-tree edge, record the number of times the
lifted relator traverses it in its chosen orientation
minus the number of times it traverses it in the
opposite orientation.

Reduction modulo $2$ shows that the first four
rows form a basis of the row space of $M$ over
$\mathbb F_2$. Hence
\[
\operatorname{rank}_{\mathbb F_2}M=4.
\]
The submatrix on rows $2,4,11,22,27,32,44$
and columns $2,3,4,5,7,8,9$, numbered from one, is
\[
\begin{pmatrix}
-1&0&-1&0&-3&1&-1\\
1&-1&-1&-2&-2&-1&0\\
0&-2&-2&-1&-3&0&0\\
-2&-1&-4&1&0&0&1\\
2&0&1&-2&-2&-2&0\\
0&0&1&-1&-5&0&-2\\
2&0&-1&-3&-3&-2&0
\end{pmatrix},
\]
and its determinant is $-8$. Therefore
$\operatorname{rank}_{\mathbb{Q}}M\geq7$.

On the other hand, the defining homomorphism
$\alpha_*:\Gamma\to\Z^2$ maps $N_2$ onto $2\Z^2$.
Consequently, $H_1(N_2;\Z)$ has free rank at least
two, so
\[
\operatorname{rank}_{\mathbb{Q}}M\leq9-2=7.
\]
Thus $\operatorname{rank}_{\mathbb{Q}}M=7$.

Let $d_1\mid\cdots\mid d_7$ be the positive
nonzero Smith invariant factors of $M$.
The rank modulo $2$ is four, so exactly four
of these factors are odd and three are even.
Hence
\[
8\mid d_1\cdots d_7.
\]
The product $d_1\cdots d_7$ is the positive
greatest common divisor of the $7\times7$ minors
of $M$. The displayed minor therefore gives
\[
d_1\cdots d_7\mid8.
\]
It follows that $d_1\cdots d_7=8$, and hence
\[
(d_1,\ldots,d_7)=(1,1,1,1,2,2,2).
\]
Since there are nine generators and seven
nonzero invariant factors, we conclude that
\[
H_1(N_2;\Z)\cong\Z^2\oplus(\Z/2)^3.
\qedhere
\]
\end{proof}

\begin{remark}\label{rmk:Magma}
$N_2$ was also found as a subgroup of the presentation of the fundamental
group of Cartwright-Steger surface in~\cite{ASY}.
Using the command \texttt{AbelianQuotient} in Magma~\cite{BCP} to the
presentation of the lattice listed at the end of the file [45] in~\cite{ASY},
the homology claimed in Lemma~\ref{lem:smaller-n2-homology} can also be verified.
\end{remark}

\begin{proof}[Proof of Theorem~{\rm \ref{thm:smaller-cyclic}}]
Fix $m\geq2$.
In $\Sigma_2\times T^2$, resolve the two positive intersections
of the sections $\Sigma_2\times\{t_1\}$ and
$\Sigma_2\times\{t_2\}$ with the fiber $\{x_0\}\times T^2$.
The resulting symplectic surface has genus $2+2+1=5$
and square $4$. Blow up one regular point of this surface
and denote its proper transform by $F'_+$. Thus
\[
B'=(\Sigma_2\times T^2)\#\CPb,\qquad
g(F'_+)=5,\qquad (F'_+)^2=3.
\]

Let us use the four Lagrangian tori and surgeries of
\cite[Section~9, Equation~(14)]{AP}, with $p=m$ and $q=1$:
\[
\begin{aligned}
&(a_1'\times c',a_1',-1),\qquad
  (b_1'\times c'',b_1',-1),\\
&(a_2'\times c',c',+1/m),\qquad
  (a_2''\times d',d',+1).
\end{aligned}
\]
Following \cite[Section~9]{AP}, let $Y_1(1/m,1)$
denote the symplectic four-manifold obtained from
$\Sigma_2\times T^2$ after these four Luttinger surgeries.

Take $t_1$ to be the standard section point of
\cite[Section~9]{AP}, and choose $t_2$ in a small disk
about $t_1$ disjoint from the projections of the surgery
tori to $T^2$. 
Choose $x_0$ away from the projections of the tori to $\Sigma_2$ and from the marked section curves.
The resolutions and blowup can then be made away from
the surgery neighborhoods and the marked curves.
Performing the same surgeries on $B'$
gives $Y_1(1/m,1)\#\CPb$. The surface $F'_+$ is
disjoint from these neighborhoods and therefore
remains embedded and symplectic after the surgeries.

In the surgered block, the first section's marked $b_2$
represents a conjugate of the generator $b_2$ in
\cite[Section~9]{AP}. The second section's marked $a_2$
represents a conjugate of $a_2$, since the two section
curves are joined by the product annulus over a path
from $t_1$ to $t_2$.

The two sections and the fiber give an orthogonal
symplectic decomposition of $H_1(F'_+;\Z)$ into
summands of ranks four, four, and two.
Let $\eta:H_1(F'_+;\Z)\to\Z$ be the homomorphism
induced by the $c$-coordinate of the torus factor.
It vanishes on the section summands and satisfies
\[
\eta([c])=1,\qquad \eta([d])=0.
\]

We choose the gluing $\psi:F'_+\to F_-$ so that
$\psi^*\ell=\eta$. This ensures that the positive
$c$-coordinate and the negative half-Albanese
coordinate agree on the summing surface after
reduction modulo $m$. These homomorphisms will
then induce an epimorphism
$\pi_1(V_m)\to\Z/m$ sending $c$ to one, which
will show that the surviving generator $c$ has
order exactly $m$. The same compatibility will
allow the gluing to lift to the corresponding
degree-$m$ covers.

Let $v$ be the oriented waist circle of the smoothing
annulus in $F_-$, and set $u=[v]\in H_1(F_-;\Z)$.
By Lemma~\ref{lem:negative-surface}, the circle $v$
is nonseparating.  
$v$ also lies in the coordinate ball where
the node is smoothed, so it is nullhomotopic in the
ambient $4$-manifold. Choosing the two blowup points
outside this ball, we obtain
\[
q_-(v)=1\in\pi_1(A)=N_2.
\]
Thus $u$ is primitive
and $\ell(u)=0$.
The homology of the normalization maps onto $2\Z^2$
in the original Albanese coordinates
\cite[proof of Proposition~6.2]{AYcurve}.
Its classes can be represented away from $v$, so
there is an $x\in u^\perp$ with $\ell(x)=1$.

Choose $u'$ with $u'\cdot u=1$. Since $x\cdot u=0$
and $\ell(x)=1$, replacing $u'$ by $u'-\ell(u')x$
preserves its intersection number with $u$ and
ensures that $\ell(u')=0$.
Since $u'\cdot u=1$, we have an orthogonal direct-sum
decomposition
\[
H_1(F_-;\Z)
=
\langle u',u\rangle\oplus\langle u',u\rangle^\perp.
\]
The intersection pairing on the second summand
is therefore unimodular.

Moreover, $s=x+(x\cdot u')u$ lies in
$\langle u',u\rangle^\perp$ and has $\ell$-value one, so the
restriction of $\ell$ to $\langle u',u\rangle^\perp$ is primitive.
By unimodularity, there is a unique class
$t\in\langle u',u\rangle^\perp$ such that
\[
\ell(y)=y\cdot t
\qquad
\text{for every }y\in\langle u',u\rangle^\perp.
\]
Since $s\cdot t=\ell(s)=1$, we obtain an orthogonal
symplectic decomposition
\[
H_1(F_-;\Z)
=
\langle u',u\rangle\oplus
\langle s,t\rangle\oplus W,
\]
where $W$ has rank six.
We have $\ell(s)=1$ and $\ell(t)=t\cdot t=0$,
while $\ell(w)=w\cdot t=0$ for every $w\in W$.
Together with $\ell(u')=\ell(u)=0$, these are
the required properties.

Next we find a primitive vector $w\in W$ with
$(q_-)_*(w)=0$, where 
\[
(q_-)_*:H_1(F_-;\Z)\longrightarrow H_1(N_2;\Z)
\]
is the homomorphism induced by $q_-$ on the abelianizations.
By Lemma~\ref{lem:smaller-n2-homology},
\[
H_1(N_2;\Z)\cong\Z^2\oplus(\Z/2)^3.
\]
Let $K$ be the kernel of the composition of
$(q_-)_*|_W$ with the projection onto the free
summand. 
Then $K$ is a direct summand of $W$
of rank at least four. 
The restriction of $(q_-)_*$
to $K$ takes values in $(\Z/2)^3$ and therefore
induces a linear map
\[
K/2K\longrightarrow(\Z/2)^3
\]
with nonzero kernel.
Lift a nonzero kernel vector using coefficients
zero or one in a basis of $K$. The resulting vector
$w$ is primitive in $K$, hence in $W$, and satisfies
$(q_-)_*(w)=0$.
Choose $w'\in W$ with $w\cdot w'=1$ and extend
$(w,w')$ to a symplectic basis of $W$.

Choose a symplectic isomorphism
$H_1(F'_+;\Z)\to H_1(F_-;\Z)$ sending the first
section's pair $([a_2],[b_2])$ to $(u',u)$,
the second section's corresponding pair to $(w,w')$,
and the fiber pair $([c],[d])$ to $(s,t)$.
Send the remaining two symplectic pairs to the
remaining pairs in $W$.
The pullback of $\ell$ is then $\eta$.

The mapping class group surjects onto the integral
symplectic group, so an orientation-preserving
diffeomorphism realizes this isomorphism.
Its image of the first section's $b_2$-curve
represents $u=[v]$. By composing on the target
with a Torelli diffeomorphism, we can arrange
that this image is isotopic to $v$ as an oriented
curve \cite[Lemma~9.4(i)]{MP}.
After matching the surface areas, Moser's theorem
gives an area-preserving representative
$\psi:F'_+\to F_-$ in the same isotopy class.
Thus $\psi^*\ell=\eta$, and the first $b_2$-curve
still maps to a curve isotopic to $v$.

The opposite normal Euler numbers give a
fiber-reversing normal-bundle gluing over $\psi$.
Form the symplectic sum and perform the four
Luttinger surgeries on the tori outside the sum
neighborhood. Denote the resulting symplectic
four-manifold by $V_m$.

Each summand contains an exceptional sphere
meeting its summing surface once.
Puncturing these spheres gives meridian disks in
the surface complements; the disk on the $F'_+$ side 
is disjoint from the Luttinger surgery neighborhoods.
Thus the complement groups are the ambient groups,
and Seifert-Van Kampen theorem gives
\[
P_m:=\pi_1(V_m)
\cong
N_2*_{\pi_1(F'_+)}\pi_1(Y_1(1/m,1)),
\]
where the map to $N_2$ is $q_-\circ\psi_*$.
Note that the blowup does not change the fundamental group
so $\pi_1(Y_1(1/m,1))\cong \pi_1(Y_1(1/m,1)\#\CPb)$. 
Since $q_-$ is surjective,
$\pi_1(Y_1(1/m,1))$ maps onto $P_m$.

By \cite[Section~9]{AP}, $\pi_1(Y_1(1/m,1))$ is generated
by $a_1,b_1,a_2,b_2,c,d$, and the following relations
hold \cite[Equation~(15)]{AP}:
\begin{gather*}
[c^{-1},b_2]=d,\qquad
[b_1^{-1},d^{-1}]=a_1,\qquad
[a_1^{-1},d]=b_1,\\
[d^{-1},b_2^{-1}]=c^m,\qquad
[a_2,c]=1.
\end{gather*}
The first section's $b_2$-curve maps to a curve
isotopic to $v$, which is trivial in $N_2$.
Hence $b_2=1$ in $P_m$, and these relations give
\[
d=a_1=b_1=1,\qquad c^m=1.
\]
Thus $P_m$ is generated by the commuting elements
$a_2,c$ and is already abelian.

The homomorphism $N_2\to P_m$ therefore factors
through $H_1(N_2;\Z)$.
The second section's $a_2$-curve has image $w$
in $H_1(F_-;\Z)$, and $(q_-)_*(w)=0$.
Its image in $P_m$ is consequently trivial.
Hence $a_2=1$, and $P_m$ is cyclic, generated by $c$,
with order dividing $m$.

To prove that $c$ has order exactly $m$, consider the
homomorphism
\[
\pi_1(\Sigma_2\times T^2)\longrightarrow\Z/m
\]
that sends $c$ to $1$ and $a_1,b_1,a_2,b_2,d$ to $0$.
Its pullback to the complement of the surgery tori
vanishes on each meridian $\mu_i$.
In the Lagrangian framings, the four filling slopes,
written additively in the first homology of the
respective boundary three-tori, are
\[
\mu_1-a_1',\qquad
\mu_2-b_1',\qquad
\mu_3+mc',\qquad
\mu_4+d'.
\]
The homomorphism sends all four filling slopes to zero
in $\Z/m$: it vanishes on the meridians and on
$a_1',b_1',d'$, while its value on $mc'$ is
$m=0$ in $\Z/m$. It therefore extends over all four
surgeries.

The homomorphism therefore extends over all four surgeries.
The map it induces on $H_1(F'_+;\Z)$ is the reduction
of $\eta$ modulo $m$. 
Recall from Section~\ref{sec:sum} that
$\vartheta_m:N_2\to\Z/m$ is the reduction of
$\vartheta$ modulo $m$. Since
$\ell=\vartheta\circ q_-$, we have
\[
\vartheta_m\circ q_-=\ell\bmod m.
\]
Since $\psi^*\ell=\eta$, the two homomorphisms
induce the same map from $\pi_1(F'_+)$ to $\Z/m$.
The universal property of the amalgamated product now gives a
homomorphism
\[
P_m=N_2*_{\pi_1(F'_+)} \pi_1(Y_1(1/m,1))\longrightarrow\Z/m.
\]
This homomorphism sends $c$ to $1$, so it is surjective.
It follows that $\pi_1(V_m)\cong\Z/m$.

The symplectic sum formulas and invariance of
Euler characteristic and signature under torus
surgery give
\[
e(V_m)=14+1+4(5-1)=31,\qquad
\sigma(V_m)=2-1=1.
\]
Therefore, we have 
\[
\chi_h(V_m)=8,\qquad c_1^2(V_m)=65.
\]

Before surgery, both summands are blowups of
aspherical manifolds. The image of any sphere
under blowdown is nullhomotopic, so its homology
class lies in the exceptional lattice.
If its square is $-1$ and its class is
$\sum_j n_jE_j$, then
\[
-1=-\sum_j n_j^2.
\]
Its class is therefore $\pm E_j$ for some $j$.
Each exceptional class meets the corresponding
summing surface once, so neither complement
contains a sphere of square $-1$.
The signatures $2$ and $-1$ exclude the
$S^2$-bundle case (with signature $0$) of Usher's theorem.
Thus the unsurgered sum is minimal
\cite[Theorem~1.1]{U}, and the four Luttinger
surgeries preserve minimality
\cite[Proposition~3.1]{HL}.

Before surgery, take the connected degree-$m$ covers
of $B'$ and $A$ defined, respectively, by the
$c$-coordinate homomorphism modulo $m$ and by
$\vartheta_m$.
The cover of $B'$ is obtained from the product cover
of $\Sigma_2\times T^2$ in the $c$ direction by
blowing up the $m$ lifts of the original blowup point.
The cover of $A$ is obtained from the cover of $X_4$
corresponding to $\ker\vartheta_m$ by blowing up the
$2m$ lifts of the two original blowup points.
The inverse images of $F'_+$ and $F_-$ are connected,
since the induced homomorphisms from their fundamental
groups to $\Z/m$ are surjective.
These lifted surfaces have genus $4m+1$ and
self-intersections $3m$ and $-3m$, respectively.
The homomorphisms from the boundary fundamental groups
to $\Z/m$ agree under the gluing and vanish on the
normal-circle fibers, so the gluing diffeomorphism
lifts. Consequently, the resulting degree-$m$ cover
of the sum before surgery is the symplectic sum of
the covered summands along the lifted surfaces.

Their blowdowns remain aspherical, and every
lifted exceptional sphere meets the lifted
summing surface once. The preceding exceptional-class
argument therefore applies to both covered summands.
Their signatures are $-m$ and $2m$, so Usher's
theorem again gives minimality before surgery.

The first three surgery tori each have one connected
preimage, covered with degree $m$ in the $c$ direction.
The meridians and the first two surgery longitudes
lift with degree one. Thus the first two filling
slopes lift to
\[
\widehat\mu_1-\widehat a_1',
\qquad
\widehat\mu_2-\widehat b_1',
\]
and both lifted surgeries have coefficient $-1$.
For the third torus, the primitive lifted longitude
$\widehat c'$ projects to $mc'$, so $\mu_3+mc'$
lifts to $\widehat\mu_3+\widehat c'$, giving
coefficient $+1$.
The homomorphism defining the cover is trivial on the
fundamental group of the fourth surgery torus.
The inverse image of this torus therefore consists of
$m$ disjoint tori, each mapping diffeomorphically onto
the original torus.
On each, the lifted
filling slope is $\widehat\mu_4+\widehat d'$
and the coefficient is $+1$.

The pulled-back framings are Lagrangian.
Hence all $m+3$ lifted surgeries are Luttinger
surgeries and preserve minimality.
Since the homomorphism defining the cover sends each
filling slope to zero in $\Z/m$, the cover extends
over all the surgeries.
Its monodromy is the isomorphism
$\pi_1(V_m)\to\Z/m$ established above.
It is therefore the universal cover
$\widetilde V_m$, with the pullback symplectic
form, and is minimal.

Since $\sigma(V_m)=1$, its unimodular intersection
form modulo torsion is odd.
Lemma~\ref{lem:cyclicparity} then implies that
$\widetilde V_m$ is nonspin.
Euler characteristic and signature multiply
under finite covers, so
\[
e(\widetilde V_m)=31m,\qquad
\sigma(\widetilde V_m)=m.
\]
Simple connectivity gives
\[
b_2^+(\widetilde V_m)=16m-1,\qquad
b_2^-(\widetilde V_m)=15m-1.
\]
Its intersection form is therefore the odd
indefinite unimodular form
\[
(16m-1)\langle1\rangle\oplus
(15m-1)\langle-1\rangle.
\]
Freedman's theorem \cite{F}, together with the
vanishing of the Kirby--Siebenmann invariant
for smooth manifolds, identifies its homeomorphism
type as
\[
(16m-1)\CP\#(15m-1)\CPb.
\]

Since $b_2^+(\widetilde V_m)>1$, Taubes's theorem
\cite{T} implies a nonzero Seiberg--Witten invariant of $\widetilde V_m$.
The standard connected sum has vanishing
Seiberg--Witten invariants, since it splits into
two summands with positive $b_2^+$.
Thus the two manifolds are not diffeomorphic.
Finally, the deck group acts freely and preserves
the pullback symplectic form.
\end{proof}

\end{document}